\documentclass[11pt,english]{article}
\usepackage[T1]{fontenc}
\usepackage[latin9]{inputenc}
\usepackage[a4paper]{geometry}
\usepackage{babel}
\usepackage{float}
\usepackage{amsthm}
\usepackage{amsmath}
\usepackage{amssymb}
\usepackage{esint}
\usepackage[authoryear]{natbib}
\usepackage{lmodern}
\usepackage{graphics}
\usepackage{microtype}
\usepackage[unicode=true,pdfusetitle,
 bookmarks=true,bookmarksnumbered=false,bookmarksopen=false,
 breaklinks=false,pdfborder={0 0 0},backref=false,colorlinks=true]
 {hyperref}

\makeatletter

\theoremstyle{plain}
\newtheorem{thm}{Theorem}
\theoremstyle{remark}
\newtheorem{rem}[thm]{Remark}

\newcommand{\includefigure}[1]{\centering\includegraphics{figures/#1}}

\def\rd{\mathrm{d}}

\def\bu{\boldsymbol{u}}
\def\bx{\boldsymbol{x}}

\def\bff{\boldsymbol{f}}
\def\bF{\boldsymbol{F}}
\def\bn{\boldsymbol{n}}
\def\be{\boldsymbol{e}}
\def\bnabla{\boldsymbol{\nabla}}
\def\bcdot{\boldsymbol{\cdot}}
\def\bwedge{\boldsymbol{\wedge}}
\def\bzero{\boldsymbol{0}}

\makeatother

\begin{document}

\title{Generalized scale-invariant solutions to the two-dimensional stationary
Navier-Stokes equations}

\author{\href{mailto:julien.guillod@unige.ch}{Julien Guillod} and \href{mailto:peter.wittwer@unige.ch}{Peter Wittwer}\\
{\small Department of Theoretical Physics,}\\
{\small University of Geneva, Switzerland}}
\maketitle
\begin{abstract}
New explicit solutions to the incompressible Navier-Stokes equations
in $\mathbb{R}^{2}\setminus\left\{ \bzero\right\} $ are determined,
which generalize the scale-invariant solutions found by Hamel. These
new solutions are invariant under a particular combination of the
scaling and rotational symmetries. They are the only solutions invariant
under this new symmetry in the same way as the Hamel solutions are
the only scale-invariant solutions. While the Hamel solutions are
parameterized by a discrete parameter $n$, the flux $\Phi$ and an
angle $\theta_{0}$, the new solutions generalize the Hamel solutions
by introducing an additional parameter $a$ which produces a rotation.
The new solutions decay like $\left|\bx\right|^{-1}$ as the Hamel
solutions, and exhibit spiral behavior. The new variety of asymptotes
induced by the existence of these solutions further emphasizes the
difficulties faced when trying to establish the asymptotic behavior
of the Navier-Stokes equations in a two-dimensional exterior domain
or in the whole plane.
\end{abstract}
\textit{\small Keywords:}{\small{} Navier-Stokes equations, Exact solutions,
Spirals, Fluid-structure interactions}\\
\textit{\small MSC class:}{\small{} 76D05, 76D03, 76U05, 35Q30, 74F10}{\small \par}

\section{Introduction}

We study a new special class of solutions to the stationary incompressible
Navier-Stokes equations in $\Omega=\mathbb{R}^{2}\setminus\left\{ \bzero\right\} $,
\begin{align}
\Delta\bu-\bnabla p & =\bu\bcdot\bnabla\bu\,, & \bnabla\bcdot\bu & =0\,, & \lim_{|\bx|\to\infty}\bu & =\bzero\,.\label{eq:ns}
\end{align}
An important parameter which labels the solutions of this system is
the flux,
\begin{equation}
\Phi=\int_{\gamma}\bu\bcdot\bn\,,\label{eq:flux}
\end{equation}
which is independent of the choice of any simple closed curve $\gamma$
encircling the origin. The equations (\ref{eq:ns}) are invariant
under two types of symmetries: the rotations around the origin $\bu(\bx)\mapsto\mathbf{R}^{-1}\bu(\mathbf{R}\bx)$,
with $\mathbf{R}\in\text{SO}(2)$ and the scaling $\bu(\bx)\mapsto\lambda\bu(\lambda\bx)$,
with $\lambda\in\mathbb{R}$. The solutions that are invariant under
these symmetries play a particular role \citep[pp. 168-173]{Wang-ExactSolutionsSteady1991}
in the asymptotic behavior of the Navier-Stokes equations, as explained
later. \citet{Sverak-LandausSolutionsNavier2011} studied in details
the scale-invariant solutions of the Navier-Stokes equations in dimension
$d\geq2$. In three dimensions, the only scale-invariant solutions
are the \citet{Landau-newexactsolution1944} solutions, which decay
like $\left|\bx\right|^{-1}$ and are labeled by a vector in $\mathbb{R}^{3}$
whose norm determines the force acting on the fluid. In two-dimensions,
\citet[§5]{Sverak-LandausSolutionsNavier2011} showed that the only
scale-invariant solutions of (\ref{eq:ns}) are the \citet[§6]{Hamel-SpiralfoermigeBewegungen1917}
solutions. The Hamel solutions are characterized by the flux $\Phi$
and a discrete parameter $n\in\mathbb{N}$, with an additional parameter
$\mu$ for $n=0$. In polar coordinates $\left(r,\theta\right)$ they
are given for $n=0$ by
\begin{equation}
\bu_{\Phi,0}=\frac{\Phi}{2\pi r}\be_{r}+\frac{\mu}{r}\be_{\theta}\,,\label{eq:hamel-n0}
\end{equation}
where $\mu\in\mathbb{R}$ is an a additional parameter, and by
\begin{equation}
\bu_{\Phi,n}=\frac{-1}{r}\varphi(\theta_{0}+\theta)\be_{r}\,,\label{eq:hamel-n}
\end{equation}
for $n\in\mathbb{N}^{*}$ and $4+\frac{\Phi}{\pi}\leq n^{2}$, where
$\varphi$ is a $\frac{2\pi}{n}$-periodic function determined by
$n$ and $\Phi$, and $\theta_{0}$ is an angle that can be chosen
arbitrarily. In view of their special form these solutions are scale-invariant,
\emph{i.e.} $\bu(\bx)=\lambda\bu(\lambda\bx)$. Moreover, it is interesting
to note that in the case $n=0$, and $\Phi\neq-4\pi$, \citet[§11]{Hamel-SpiralfoermigeBewegungen1917}
found one more free parameter $A\in\mathbb{R}$ since
\begin{equation}
\bu_{\Phi,0,A}=\frac{\Phi}{2\pi r}\be_{r}+\left(\frac{\mu}{r}+Ar^{1+\frac{\Phi}{2\pi}}\right)\be_{\theta}\label{eq:hamel-n0-A}
\end{equation}
is an exact solution of (\ref{eq:ns}) provided $\Phi<-2\pi$. This
solution in not scale-invariant and is bounded by $r^{-1}$ at infinity
only for $\Phi\leq-4\pi$.

In what follows we look for solutions invariant under combinations
of the scaling and rotational symmetries. We say that a solution $\bu$
of the Navier-Stokes equations (\ref{eq:ns}) is scale-invariant up
to a rotation if there exists a rotation matrix $\mathbf{R}_{\lambda}\in\text{SO}(2)$
of continuously differentiable angle $R(\lambda)$ such that
\begin{equation}
\lambda\boldsymbol{u}(\lambda\bx)=\mathbf{R}_{\lambda}^{-1}\boldsymbol{u}(\mathbf{R}_{\lambda}\bx)\,,\label{eq:symmetry}
\end{equation}
for all $\lambda>0$. The scale-invariant solutions corresponds to
the special case $R(\lambda)\equiv0$. The aim of this paper is to
determine all solutions of (\ref{eq:ns}) that are scale-invariant
up to a rotation, with $R(\lambda)$ a continuously differentiable
function and discuss their implications.

Our main result is the following:
\begin{thm}
For all $n\in\mathbb{N}^{*}$, $\Phi\in\mathbb{R}$ and $a\in\mathbb{R}$
satisfying
\[
\frac{4+\frac{\Phi}{\pi}}{1+a^{2}}\leq n^{2}\,,
\]
there exists a $\frac{2\pi}{n}$-periodic function $\varphi$ depending
on $n$, $\Phi$, and $a$, such that for any $\theta_{0}\in\mathbb{R}$,
\begin{equation}
\bu_{n,\Phi,a}=\frac{1}{r}\left[-\varphi(\theta_{0}+\theta+a\log r)\be_{r}+a\left(\varphi(\theta_{0}+\theta+a\log r)-4\right)\be_{\theta}\right]\,,\label{eq:u-new}
\end{equation}
and the associated pressure (\ref{eq:p}) satisfy the Navier-Stokes
equations (\ref{eq:ns}). These solutions are invariant under the
symmetry (\ref{eq:symmetry}) with $R(\lambda)=-a\log\lambda$ and
have flux $\Phi$. Moreover any solution of the Navier-Stokes equations
(\ref{eq:ns}) which is invariant under the symmetry (\ref{eq:symmetry})
for some continuously differentiable rotation $R(\lambda)$ is equal
either to one of the exact solution $\bu_{n,\Phi,a}$ for an angle
$\theta_{0}$, or to a Hamel solution $\bu_{\Phi,0}$ defined by (\ref{eq:hamel-n0})
with $n=0$ for some $\mu\in\mathbb{R}$.\end{thm}
\begin{rem}
The ansatz for spiral solutions made by \citet[§9]{Hamel-SpiralfoermigeBewegungen1917}
does not allow solutions in the plane with streamlines that are logarithmic
spirals. The solutions with logarithmic spirals that he found are
only possible between two walls of logarithmic shape. The solutions
presented here essentially correspond to the intuition of Hamel to
look for non-harmonic function in the plane having streamlines that
are spirals.
\end{rem}

\begin{rem}
The expression (\ref{eq:u-new}) is a solution of the Navier-Stokes
equations in $\mathbb{R}^{2}\setminus\left\{ \bzero\right\} $, but
due to the behavior near the origin like $r^{-1}$, the non-linear
term $\bu\bcdot\bnabla\bu$, even when written as $\bnabla\bcdot\left(\bu\otimes\bu\right)$,
has no immediate distributional meaning in $\mathbb{R}^{2}$. This
is in contrast to the three-dimensional case where the non-linear
term of a scale-invariant solution is a distribution even if $\bu$
diverges likes $r^{-1}$ at the origin. One can nevertheless always
construct a solution to the Navier-Stokes equations in $\mathbb{R}^{2}$
by truncating one of the exact solutions near the origin, and defining
the source term by the truncation error. The force and the torque
of a solution $\bu$ are given for any curve $\gamma$ encircling
the origin, by
\begin{align*}
\bF & =\int_{\gamma}\mathbf{T}\bn\,, & M & =\int_{\gamma}\mathbf{x}\bwedge\mathbf{T}\bn\,,
\end{align*}
where $\mathbf{T}$ is the stress tensor including the convective
part, $\mathbf{T}=\bu\otimes\bu+p-\bnabla\bu-\left(\bnabla\bu\right)^{T}$.
By taking for $\gamma$ a circle whose radius goes to infinity, the
force is zero, $\bF=\bzero$. By taking for simplicity the circle
of radius one, the torque is 
\[
M=a\left(16\pi+6\Phi+\int_{-\pi}^{+\pi}\varphi^{2}(\theta)\,\rd\theta\right)\,.
\]

\end{rem}
The study of scale-invariant solutions has proven to be of great importance,
in particular for the determination of the asymptotic behavior of
the stationary Navier-Stokes equations in two or three dimensions.
The stationary and incompressible Navier-Stokes equations in the exterior
domain $\Omega=\mathbb{R}^{2}\setminus B$ of a compact, connected
set $B$ are\begin{subequations}
\begin{equation}
\begin{aligned}\Delta\bu-\bnabla p & =\bu\bcdot\bnabla\bu\,, & \bnabla\bcdot\bu & =0\,,\\
\left.\bu\right|_{\partial B} & =\bu^{*}\,, & \lim_{|\bx|\to\infty}\bu & =\bzero\,,
\end{aligned}
\label{eq:ns-0}
\end{equation}
where $\bu^{*}$ is any smooth boundary condition with no net flux,
\begin{equation}
\int_{\partial B}\bu^{*}\bcdot\bn=0\,.\label{eq:no-flux}
\end{equation}
\label{eq:ns-noforce}\end{subequations} Problem (\ref{eq:ns-noforce})
is closely related to the one of the incompressible Navier-Stokes
equations in $\mathbb{R}^{2}$,
\begin{align}
\Delta\bu-\bnabla p & -\bu\bcdot\bnabla\bu=\bff\,, & \bnabla\bcdot\bu & =0\,, & \lim_{|\bx|\to\infty}\bu & =\bzero\,,\label{eq:ns-force}
\end{align}
where $\bff$ is a smooth function of compact support. We remark,
that the problems (\ref{eq:ns-noforce}) and (\ref{eq:ns-force})
are very similar on a formal level: any solution of (\ref{eq:ns-noforce})
defines a solution of (\ref{eq:ns-force}) on the exterior of the
support of $\bff$, and conversely any solution of (\ref{eq:ns-noforce})
can be truncated in order to obtain a solution of (\ref{eq:ns-force}).
In three dimensions, \citet{Nazarov.Pileckas-AsymptoticofSolutions1999,Nazarov-steady2000}
proved that the asymptotic behavior of solutions of (\ref{eq:ns-noforce})
is a scale-invariant solution. Then \citet{Korolev.Sverak-largedistanceasymptotics2011}
simplified the proof by showing directly that, in this case, the Landau
solution is the correct asymptotic behavior of any solution bounded
by $\left(1+\left|\bx\right|\right)^{-1}$. In two dimensions, existence
of solutions to (\ref{eq:ns-noforce}) or (\ref{eq:ns-force}) are
not known in general \citep{Galdi-StationaryNavier-Stokesproblem2004,Guillod-Asymptoticbehaviour2013},
even for small data. The difference between two and three dimensions
is essentially that in three dimensions the compatibility condition
of the Stokes approximation to decay faster than $r^{-1}$ at infinity
corresponds to the force and can be lifted by the Landau solutions
which are exact solution of (\ref{eq:ns-force}) with $\bff(\bx)=\boldsymbol{b}\delta(\bx)$
where $\boldsymbol{b}\in\mathbb{R}^{3}$ is the net force and $\delta$
is the Dirac distribution. In two dimensions and in the case where
$\bff$ has non-zero mean, \citet{Guillod-Asymptoticbehaviour2013}
showed by physical arguments and detailed numerical verification that
the velocity has to decay like $r^{-1/3}$ at infinity. In the case
where $\bff$ has zero mean, one would guess by analogy with the three-dimensional
case that the asymptotic behavior should be a scale-invariant solution.
However, \citet[§5]{Sverak-LandausSolutionsNavier2011} shows that
one can not prove this by using perturbation techniques based on the
Stokes approximation, and even together with the newly discovered
solutions, we do not appear to be able to parameterize the general
asymptotic behavior in the case where $\bff$ has zero mean. The intuitive
reason for this is the fact that the Stokes approximation has two
compatibility conditions if we require the solution to decay faster
than $r^{-1}$ at infinity: one of them might be lifted by adjusting
the parameter $a$ of the new solutions, but we do not have sufficient
parameters to lift also the other compatibility condition. We believe
that the newly discovered solutions are a special case of a more general
family of solutions, yet to be discovered, with one more parameter,
corresponding to the general asymptotic behavior in the case where
$\bff$ has zero mean. 

The paper in organized as follow. We first prove that the solutions
which are scale-invariant up to a rotation are given explicitly in
term of a $2\pi$-periodic function $\varphi$ satisfying an ordinary
differential equation, and then we solve this differential equation
by using elliptic functions. Finally we represent the solutions graphically,
analyze the solutions having small amplitude, and discuss the implications
for the solutions of Navier-Stokes equations.

\section{Reduction to an ordinary differential equation}

We consider a solution $\bu$ which is scale-invariant up to a rotation
as defined in (\ref{eq:symmetry}). In polar coordinates $\left(r,\theta\right)$
this symmetry is more easily expressed,
\begin{align*}
\lambda u_{r}(\lambda r,\theta) & =u_{r}(r,\theta-R(\lambda))\,, & \lambda u_{\theta}(\lambda r,\theta) & =u_{\theta}(r,\theta-R(\lambda))\,.
\end{align*}
Therefore, by setting $\lambda=r^{-1}$, $u_{r}$ and $u_{\theta}$
are characterized by their values on $S^{1}$,
\begin{align*}
u_{r}(r,\theta) & =\frac{1}{r}\varphi_{r}(\theta+R(r^{-1}))\,, & u_{\theta}(r,\theta) & =\frac{1}{r}\varphi_{\theta}(\theta+R(r^{-1}))\,,
\end{align*}
where $\varphi_{i}(\theta)=u_{i}(1,\theta)$ for $i\in\left\{ r,\theta\right\} $.
The divergence of the vector field $\bu=u_{r}\be_{r}+u_{\theta}\be_{\theta}$
is
\[
\bnabla\bcdot\bu=\frac{1}{r^{2}}\left[\varphi_{\theta}^{\prime}(z)-r^{-1}R^{\prime}(r^{-1})\varphi_{r}^{\prime}(z)\right]\,,
\]
where $z=\theta+R(r^{-1})$. The requirement of $\bu$ to be divergence
free therefore implies that
\begin{align*}
R(\lambda) & =\theta_{0}-a\log\lambda\,, & \varphi_{\theta}(z) & =\mu-a\varphi_{r}(z)\,,
\end{align*}
where $a$, $\theta_{0}$ and $\mu$ are real constants. Consequently,
a divergence free vector field satisfies the symmetry (\ref{eq:symmetry}),
if and only if $R(\lambda)=\theta_{0}-a\log\lambda$, and if it has
the form
\begin{align}
\bu & =\frac{1}{r}\left[-\varphi(z)\be_{r}+\left(\mu+a\varphi(z)\right)\be_{\theta}\right]\,, & z & =\theta_{0}+\theta+a\log r\,,\label{eq:u}
\end{align}
where $\varphi$ is a $2\pi$-periodic function. The corresponding
stream function $\psi$, defined such that $\bu=\bnabla\bwedge\psi$,
is
\[
\psi(r,\theta)=\mu\log r+\Gamma(z)\,,
\]
where $\Gamma$ is an antiderivative of $\varphi$.

We now determine the ordinary differential equation which $\varphi$
has to satisfy in order for $\bu$ to be an exact solution of (\ref{eq:ns}).
The vorticity is
\[
\omega=\frac{1+a^{2}}{r^{2}}\varphi^{\prime}(z)\,,
\]
and the vorticity equation
\[
\Delta\omega=\bu\bcdot\bnabla\omega
\]
becomes, after an explicit integration, the following ordinary differential
equation
\begin{equation}
\left(1+a^{2}\right)\varphi^{\prime\prime}(z)-\left(\mu+4a\right)\varphi^{\prime}(z)+4\varphi(z)=\varphi(z)^{2}-C\,,\label{eq:ode}
\end{equation}
where $C\in\mathbb{R}$ is a constant related to certain averages
of $\varphi$. By integrating the Navier-Stokes equations we can construct
the pressure,
\begin{equation}
p=\frac{1}{r^{2}}\left[a\left(1+a^{2}\right)\varphi^{\prime}(z)-\left(2\left(1+a^{2}\right)+a\mu\right)\varphi(z)\right]-\frac{1}{2r^{2}}\left[\mu^{2}+C\left(1+a^{2}\right)\right]\,.\label{eq:p}
\end{equation}
This shows that the only solutions of (\ref{eq:ns}) which are scale-invariant
up to a rotation are given by (\ref{eq:u}) and (\ref{eq:p}) where
$\varphi$ is a $2\pi$-periodic function satisfying (\ref{eq:ode}).
The differential equation (\ref{eq:ode}) is analog to the one describing
the motion of a particle in a potential undergoing friction, and in
order to obtain periodic solutions, the damping term has to vanish,
\emph{i.e.} $\mu+4a=0$. So we finally end up with a differential
equation with two parameters: $a$ and $C$. In the next section we
find the periodic solutions of this differential equation.

\section{Resolution of the ordinary differential equation}

The ordinary differential equation (\ref{eq:ode}) is clearly invariant
under the translation $z\mapsto z+\theta_{0}$, so we do not keep
track of this trivial symmetry and fix the origin later on in a convenient
way. The trivial solutions where $\varphi$ is constant are not included
in this analysis, since they correspond to the Hamel solutions (\ref{eq:hamel-n0}).
As explained above, in order to obtain periodic solutions we have
to take $\mu+4a=0$, and therefore the ordinary differential equation
(\ref{eq:ode}) can be written as the differential equation describing
a free particle in a potential,
\begin{align*}
\varphi^{\prime\prime} & =-V^{\prime}(\varphi)\,, & V(\varphi) & =\frac{1}{1+a^{2}}\left[C\varphi+2\varphi^{2}-\frac{\varphi^{3}}{3}\right]\,.
\end{align*}
Finding the solutions of such an equation is rather standard \citep[Theorem 2]{Rosenhead-SteadyTwo-DimensionalRadial1940,Sverak-LandausSolutionsNavier2011}.
The energy is conserved so
\begin{align*}
E & =\frac{1}{2}\left(\varphi^{\prime}\right)^{2}+V(\varphi)\,, & \varphi^{\prime} & =\pm\sqrt{2E-2V(\varphi)}\,.
\end{align*}
Since we look for non-trivial periodic solutions, the potential has
to have a minimum, so $C>-4$, and the energy has to be between the
maximum and the minimum admissible values,
\[
\frac{2}{3}\left(\sqrt{C+4}-2\right)\left(\sqrt{C+4}-2-C\right)<\left(1+a^{2}\right)E<\frac{2}{3}\left(\sqrt{C+4}+2\right)\left(\sqrt{C+4}+2+C\right)\,.
\]
These two conditions imply that the polynomial $2E-2V(\varphi)$ has
three distinct real roots, $\varphi_{1}<\varphi_{2}<\varphi_{3}$,
and by Vieta's formulas,
\begin{equation}
\varphi_{1}+\varphi_{2}+\varphi_{3}=6\,.\label{eq:viete}
\end{equation}
Therefore,
\[
2E-2V(\varphi)=\frac{2}{3\left(1+a^{2}\right)}\left(\varphi-\varphi_{1}\right)\left(\varphi-\varphi_{2}\right)\left(\varphi-\varphi_{3}\right)\,,
\]
and the solution is given in term of the incomplete elliptic function
of the first kind $F$,
\[
z(\varphi)=\sqrt{\frac{3}{2}}\sqrt{1+a^{2}}\int_{\varphi_{1}}^{\varphi}\frac{\rd\varphi}{\sqrt{\left(\varphi-\varphi_{1}\right)\left(\varphi-\varphi_{2}\right)\left(\varphi-\varphi_{3}\right)}}=\sqrt{6}\frac{\sqrt{1+a^{2}}}{\sqrt{\varphi_{3}-\varphi_{1}}}F\left(\sqrt{\frac{\varphi-\varphi_{1}}{\varphi_{2}-\varphi_{1}}};\alpha\right)\,,
\]
where
\[
\alpha=\sqrt{\frac{\varphi_{2}-\varphi_{1}}{\varphi_{3}-\varphi_{1}}}\,.
\]
We take the following convention for the elliptic integral $F$ \citep[§8.1]{Gradshteyn-TableofIntegrals2007},
\[
F(x;\alpha)=\int_{0}^{x}\frac{\rd t}{\sqrt{\left(1-t^{2}\right)\left(1-\alpha^{2}t^{2}\right)}}=\int_{0}^{\arcsin x}\frac{\rd\theta}{\sqrt{1-\alpha^{2}\sin^{2}\theta}}\,.
\]
The function $\varphi$ is $2\pi$-periodic if there exists $n\in\mathbb{N}^{*}$
such that
\[
z(\varphi_{2})=\frac{\pi}{n}\,,
\]
\emph{i.e.,} explicitly,
\begin{equation}
\frac{2\sqrt{1+a^{2}}K(\alpha)}{\sqrt{\varphi_{3}-\varphi_{1}}}=\sqrt{\frac{2}{3}}\frac{\pi}{n}\,,\label{eq:cond-per}
\end{equation}
where $K$ is the complete elliptic function of the first kind. The
flux is given by
\[
\Phi=\int_{\partial\Omega}\bu\bcdot\bn=\int_{\partial B(\bzero,1)}\bu\bcdot\be_{r}=-\int_{-\pi}^{+\pi}\varphi(z)\,\rd z=-2n\int_{0}^{\pi/n}\varphi(z)\,\rd z\,,
\]
and, explicitly, by using the complete elliptic function of the second
kind $E$,
\begin{equation}
\Phi=-n\sqrt{6}\int_{\varphi_{1}}^{\varphi_{2}}\frac{\sqrt{1+a^{2}}\,\varphi\,\rd\varphi}{\sqrt{\left(\varphi-\varphi_{1}\right)\left(\varphi-\varphi_{2}\right)\left(\varphi-\varphi_{3}\right)}}=-2\sqrt{6}n\frac{\sqrt{1+a^{2}}}{\sqrt{\varphi_{3}-\varphi_{1}}}\left[\varphi_{3}K(\alpha)-\left(\varphi_{3}-\varphi_{1}\right)E(\alpha)\right]\,.\label{eq:cond-flux}
\end{equation}
The conditions (\ref{eq:viete}), (\ref{eq:cond-per}) and (\ref{eq:cond-flux})
reduce to
\begin{align}
H(\alpha) & =\frac{1}{n^{2}\left(1+a^{2}\right)}\left(\pi^{2}+\frac{\pi\Phi}{4}\right)\,, & H(\alpha) & =\left[\left(\alpha^{2}-2\right)K(\alpha)+3E(\alpha)\right]K(\alpha)\,.\label{eq:H}
\end{align}
Since the function $H$ is monotonic for $\alpha>0$ and its image
is $\left(-\infty,\frac{\pi^{2}}{4}\right]$, this equation has a
unique solution $\alpha_{n}>0$ for each $n\in\mathbb{N}^{*}$ satisfying
\begin{equation}
\frac{4+\frac{\Phi}{\pi}}{1+a^{2}}\leq n^{2}\,.\label{eq:cond}
\end{equation}
Since the equation determining $\alpha$ is continuous with respect
to $a$ and $\Phi$, the solution $\varphi$ depends continuously
on $a$ and $\Phi$ inside the region defined by (\ref{eq:cond}).

\section{Discussion of solutions}

\begin{figure}[H]
\includefigure{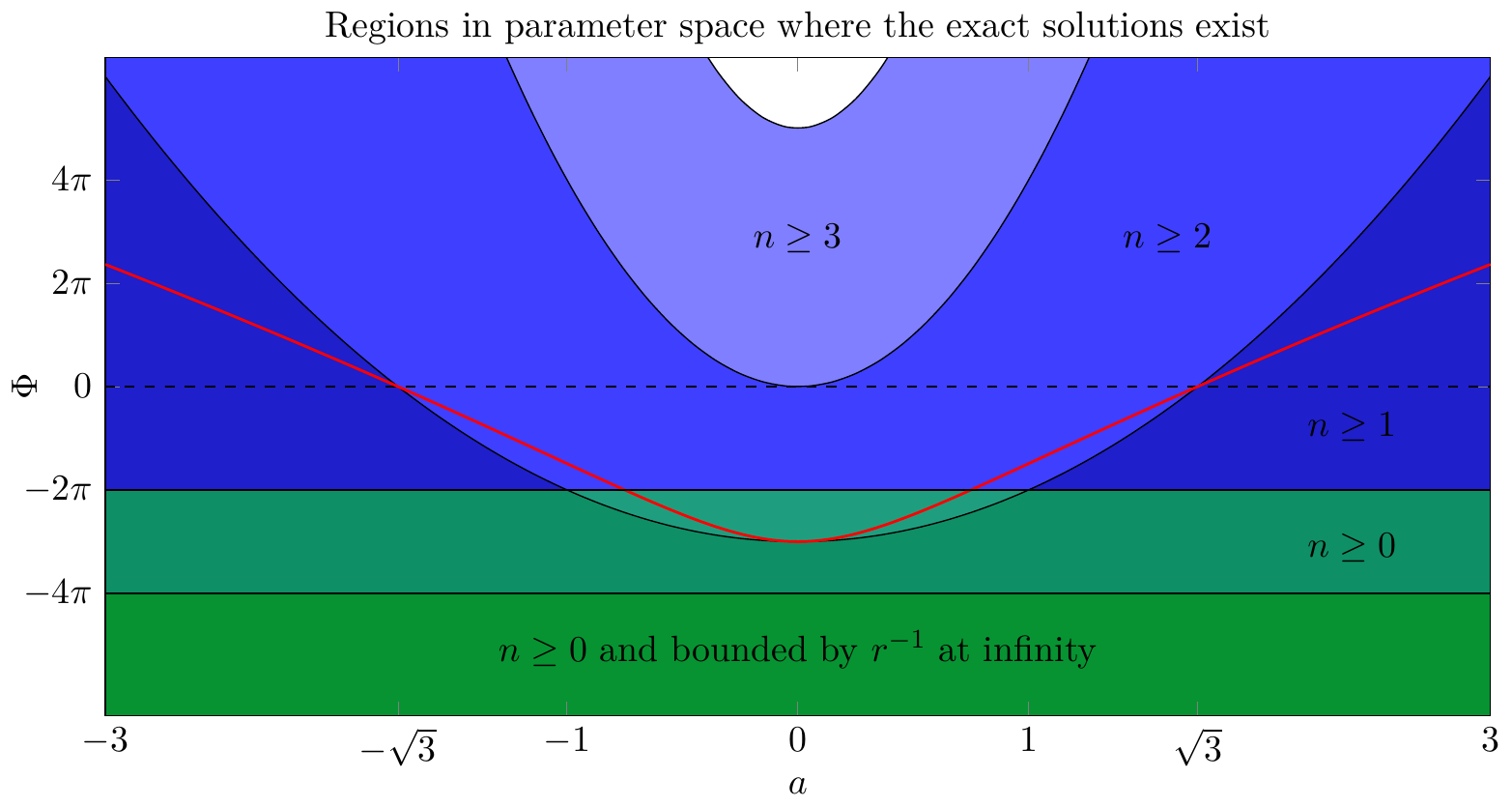}

\caption{\label{fig:diagram}Regions in the $\left(a,\Phi\right)$-plane where
the exact solutions $\bu_{\Phi,n,a}$ and $\bu_{\phi,0,A}$ exist.
For $n\geq1$, the exact solutions $\bu_{\Phi,n,a}$ exist in the
region below the parabolas filled in blue. For $n=0$, the solution
$\bu_{\phi,0,A}$ exists for $\Phi<-2\pi$ and decays like $r^{-1}$
if $\Phi\leq-4\pi$; these regions are colored in green. The red curve
represents the critical line above which the linearization of the
Navier-Stokes equations around the harmonic function (\ref{eq:hamel-n0})
with $\mu=-4a$ decays to infinity faster than $r^{-1}$.}
\end{figure}

For $n\in\mathbb{N}^{*}$, the exact solution $\bu_{\Phi,n,a}$ exists
provided condition (\ref{eq:cond}) is satisfied. The corresponding
region in the plane $\left(a,\Phi\right)$ is represented in figure~\ref{fig:diagram}.
Moreover, for $n=0$, the Hamel solution $\bu_{\Phi,0,A}$ defined
by (\ref{eq:hamel-n0-A}) with $A\neq0$ exists for $\Phi<-2\pi$
and decays like $r^{-1}$ at infinity if $\Phi\leq-4\pi$, and less
rapidly if $-4\pi<\Phi<-2\pi$. The linearization of the Navier-Stokes
equations around the harmonic function (\ref{eq:hamel-n0}) with $\mu=-4a$,
can be solved exactly by the use of a Fourier series \citep{Hillairet-mu2013}.
The Fourier modes which are not zero decay faster than $r^{-1}$ provided
a condition on $\Phi$ and $\mu$ holds. This condition is represented
in figure~\ref{fig:diagram} by a red curve, so that for all values
above the curve the Fourier modes decay faster than $r^{-1}$ at infinity.
The zero-flux case, was treated by \citet{Hillairet-mu2013}, and
they found that provided $\left|a\right|>\sqrt{3}$, the Navier-Stokes
equations in the exterior of a unit disk with a Dirichlet boundary
condition sufficiently close to $-4a\be_{\theta}$ admits a solution
whose asymptote in given by $-4\tilde{a}\be_{\theta}/r$ where $\tilde{a}$
is close to $a$.

\begin{figure}[H]
\includefigure{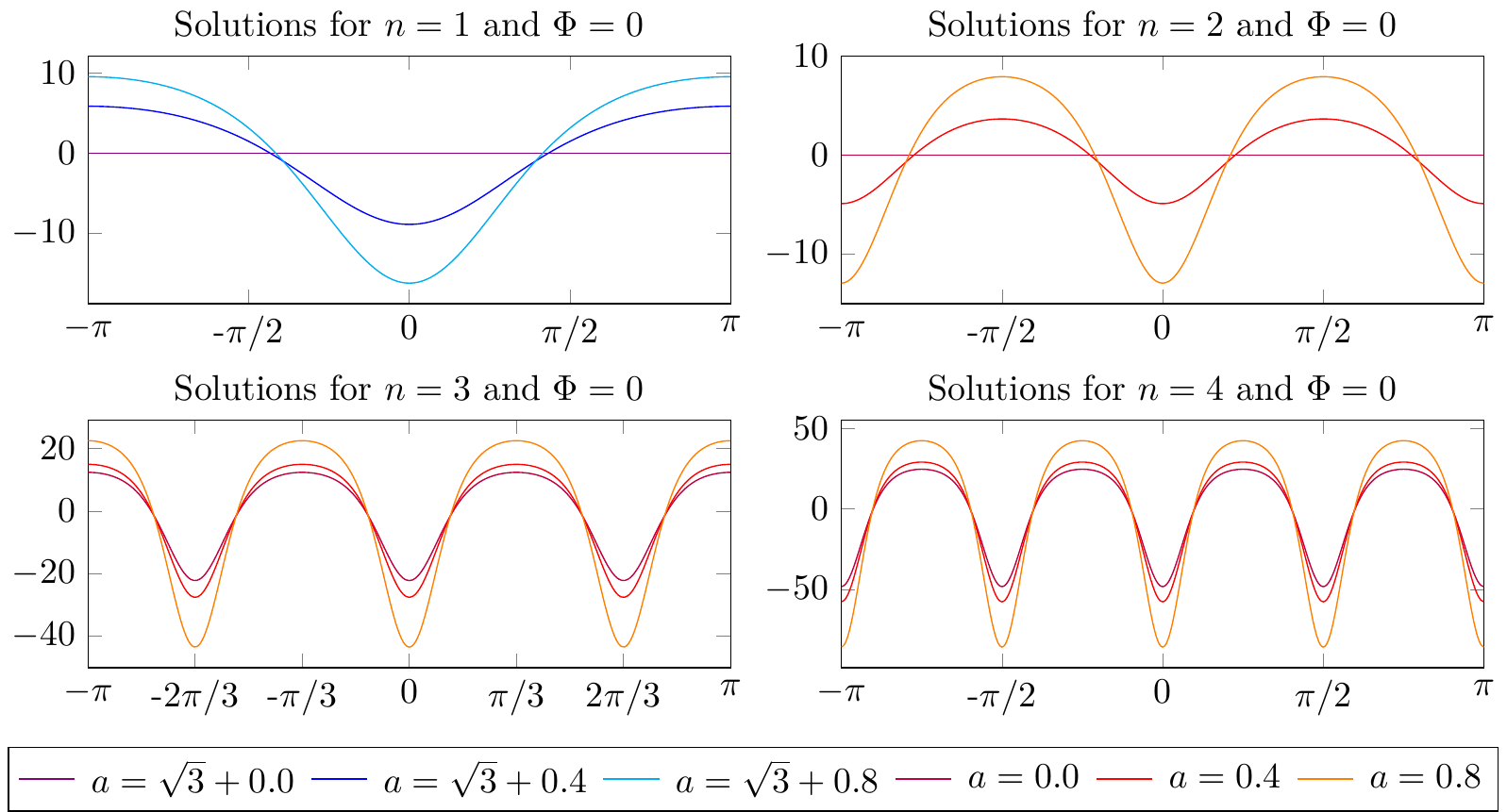}

\caption{\label{fig:ode}Periodic solutions $\varphi$ of the differential
equation (\ref{eq:ode}) for $\Phi=0$ and $n\in\left\{ 1,2,3,4\right\} $.
For $n=1$, the solution start to exists for $\left|a\right|>\sqrt{3}$,
so that is why in this case the values of $a$ start at $\sqrt{3}$.
As shown above the solutions are $\frac{2\pi}{n}$-periodic.}
\end{figure}

\begin{figure}[H]
\includefigure{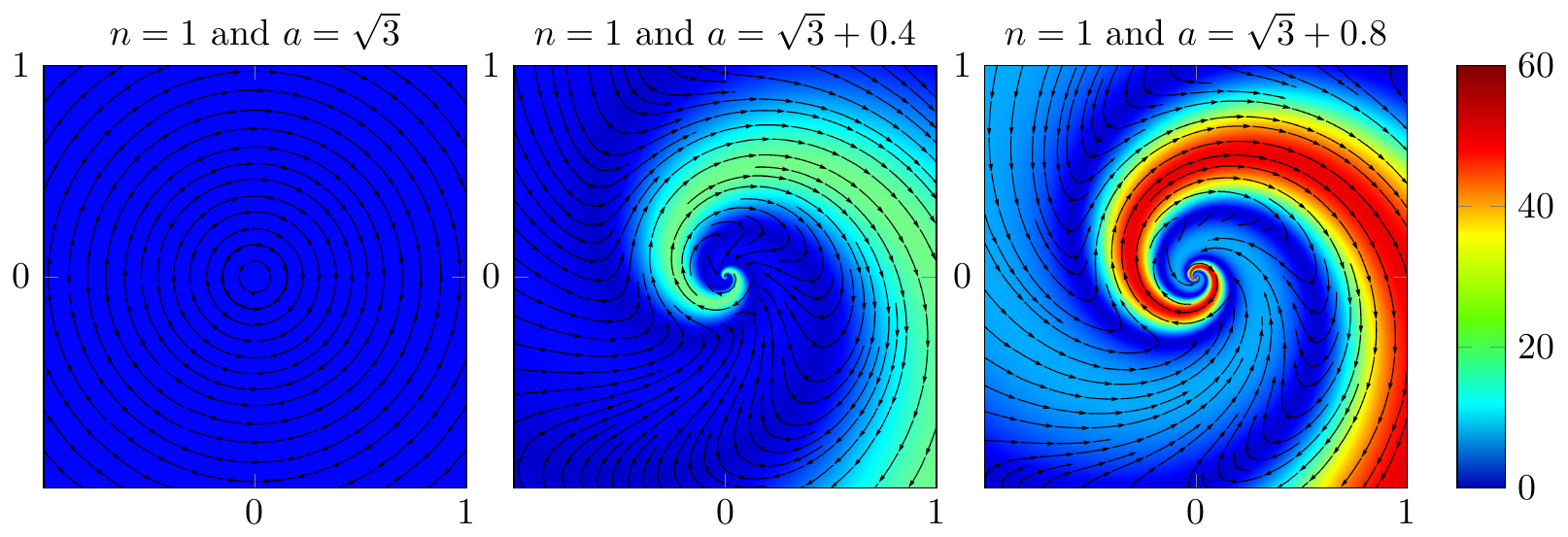}

\includefigure{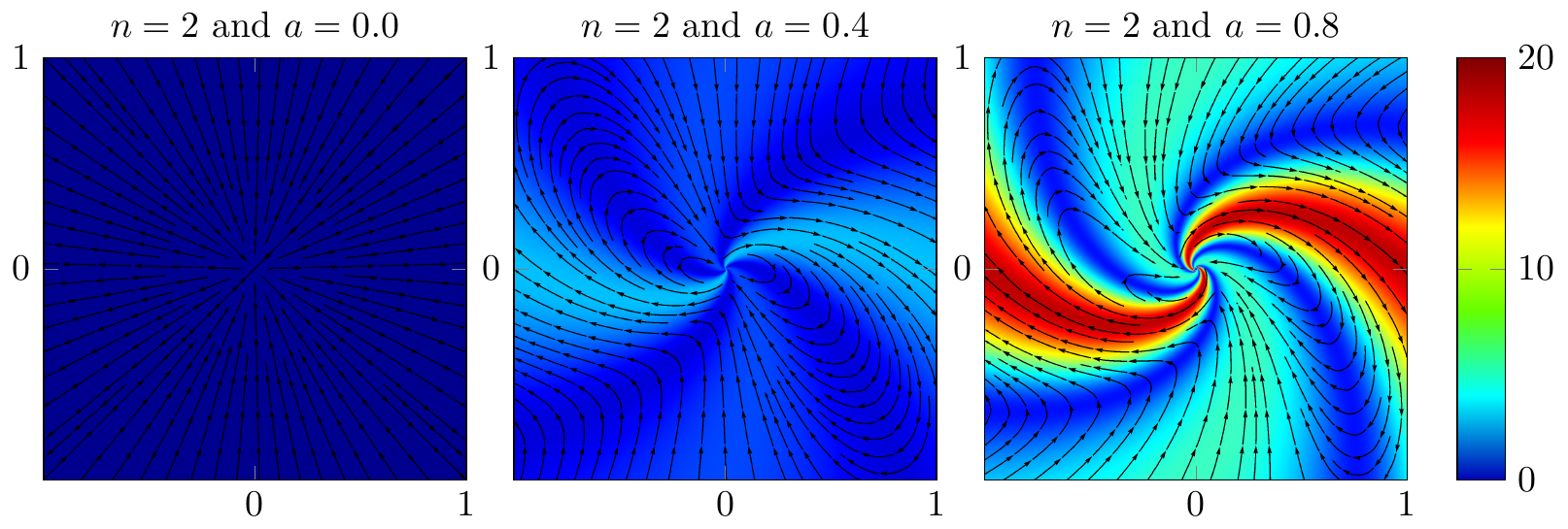}

\includefigure{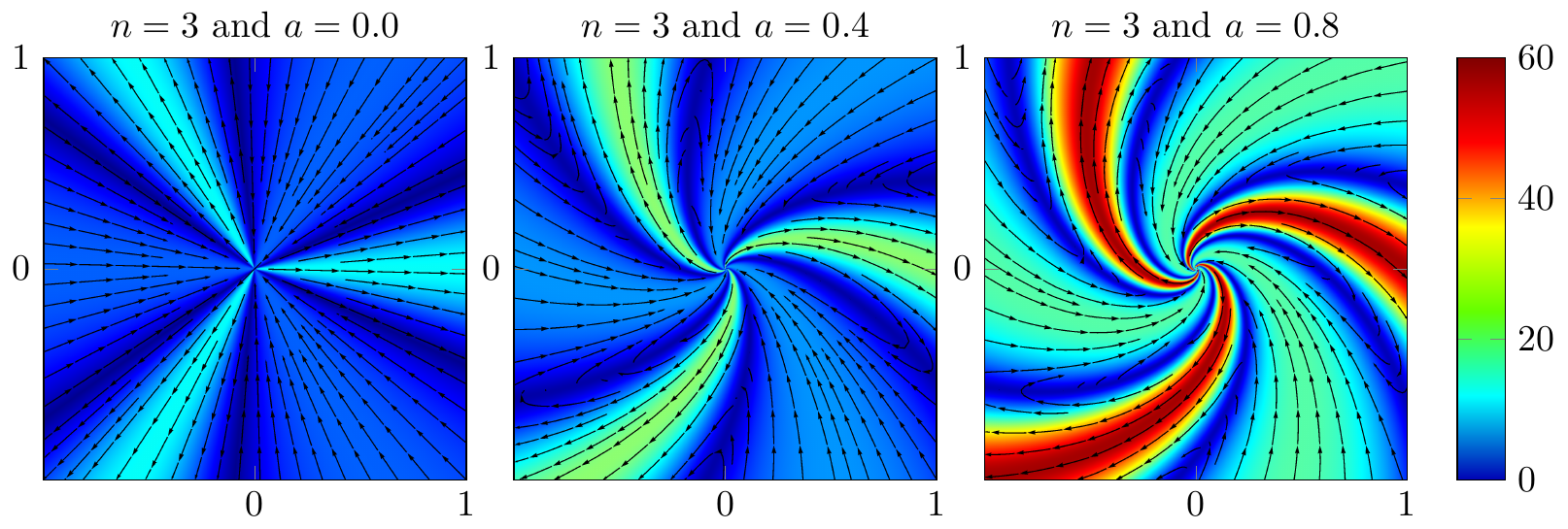}

\includefigure{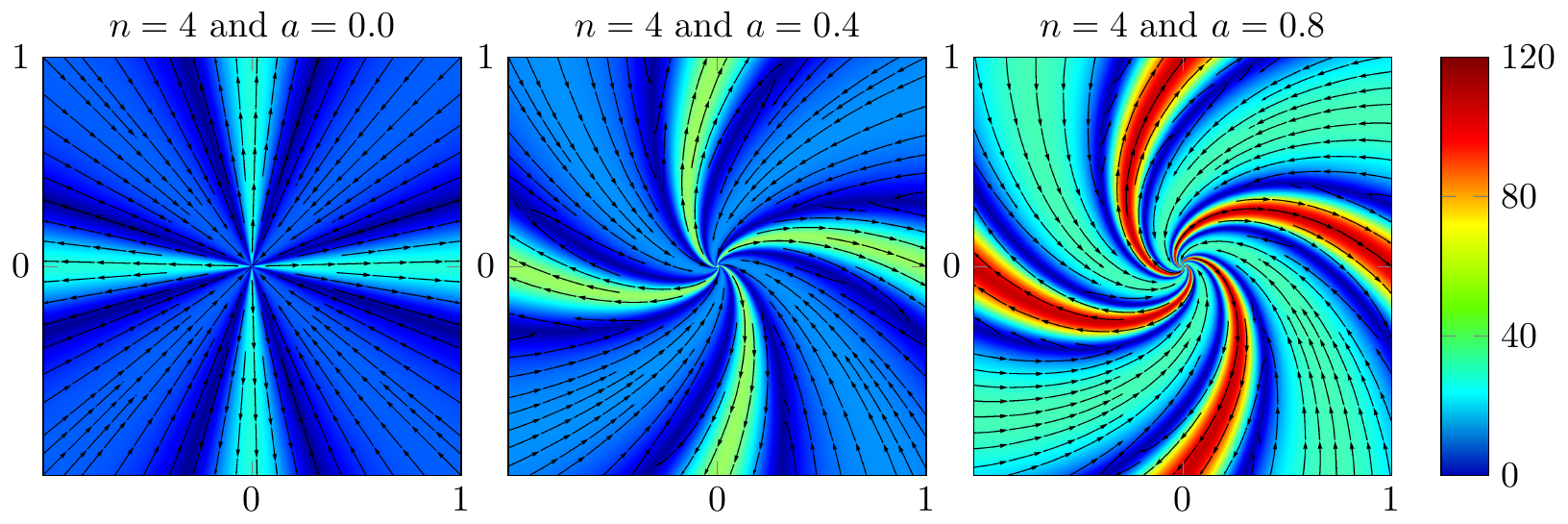}

\caption{\label{fig:plot}Representation of the velocity vector field (\ref{eq:u-new})
in case of zero flux, $\Phi=0$, for different values of $n$ and
$a$. The black lines represent the streamlines and the color the
strength of the field $r\left|\bu\right|$. By increasing the value
of $\left|a\right|$, the $n$ branches where the velocity is high
rotate more rapidly.}
\end{figure}

For a given $n$, the solution $\varphi$ has $n$ maxima and $n$
minima (figure~\ref{fig:ode}) and the parameter $a$ has the effect
of rotating the branches corresponding to these maxima or minima,
as shown in figure~\ref{fig:plot}. Small solutions of (\ref{eq:ns})
are of particular interest because even in this case we don't know
the existence of a solution in general. Certain large solutions of
the Navier-Stokes equations are known to exhibit exotic behavior at
infinity, like for example the Hamel solutions which have arbitrary
slow decay to infinity and violate uniqueness \citep[XII.2]{Galdi-IntroductiontoMathematical2011}.
Small solutions of the ordinary differential equation are given by
$\alpha$ small, and to discuss these solutions we develop $H(\alpha)$
in a series,
\[
H(\alpha)=\frac{\pi^{2}}{4}\left(1-\frac{3}{32}\alpha^{4}\right)+O(\alpha^{6})\,.
\]
In view of (\ref{eq:H}), small solutions with $\Phi=0$ are only
possible for $n\in\left\{ 1,2\right\} $. For $n=1$ and $\Phi=0$,
the solutions are defined for $\left|a\right|\geq\sqrt{3}$, and we
take for example $a=\sqrt{3}+\varepsilon$, with $\varepsilon>0$.
We find by a series expansion, that 
\[
\alpha=\left(\frac{256}{3}\right)^{1/8}\varepsilon^{1/4}+O(\varepsilon^{1/2})\,,
\]
and
\begin{align*}
\varphi_{1} & =-4\,3^{3/4}\varepsilon^{1/2}+O(\varepsilon)\,, & \varphi_{2} & =4\,3^{3/4}\varepsilon^{1/2}+O(\varepsilon)\,, & \varphi_{3} & =6+O(\varepsilon)\,,
\end{align*}
and the solution satisfies
\[
\varphi(z)=-4\,3^{3/4}\varepsilon^{1/2}\cos(z)+O(\varepsilon)\,.
\]
Since in this case $\left|a\right|$ has to be large, we note that
this does not produce a solution having a small velocity field, and
the torque is also large, $M=16\pi\sqrt{3}+O(\varepsilon)$. Moreover,
since $\mu+4a=0$, this corresponds to $\left|\mu\right|\geq\sqrt{48}$.
This specific value is interesting since this is exactly the criterion
found by \citet{Hillairet-mu2013} to obtain a solution of the Navier-Stokes
equations (\ref{eq:ns-noforce}) having the asymptote $\mu\be_{\theta}/r$
for the velocity. We note that this is not in contradiction with the
existence of the exact solutions found here, because the boundary
condition given by the evaluation of the exact solution for $a=\sqrt{3}+\varepsilon$
is too big for the theorem of \citet{Hillairet-mu2013} to apply.

For $n=2$, we can take $a=\varepsilon$, so
\[
\alpha=\left(\frac{32}{3}\right)^{1/4}\varepsilon^{1/2}+O(\varepsilon)\,,
\]
and
\begin{align*}
\varphi_{1} & =-4\sqrt{6}\,\varepsilon+O(\varepsilon^{2})\,, & \varphi_{2} & =4\sqrt{6}\,\varepsilon+O(\varepsilon^{2})\,, & \varphi_{3} & =6+O(\varepsilon^{2})\,,
\end{align*}
and the solution satisfies
\[
\varphi(z)=-4\sqrt{6}\,\varepsilon\cos(2z)+O(\varepsilon^{2})\,.
\]
The torque of this solution is given by $M=16\pi\varepsilon+O(\varepsilon^{2})$.

\begin{figure}[H]
\includefigure{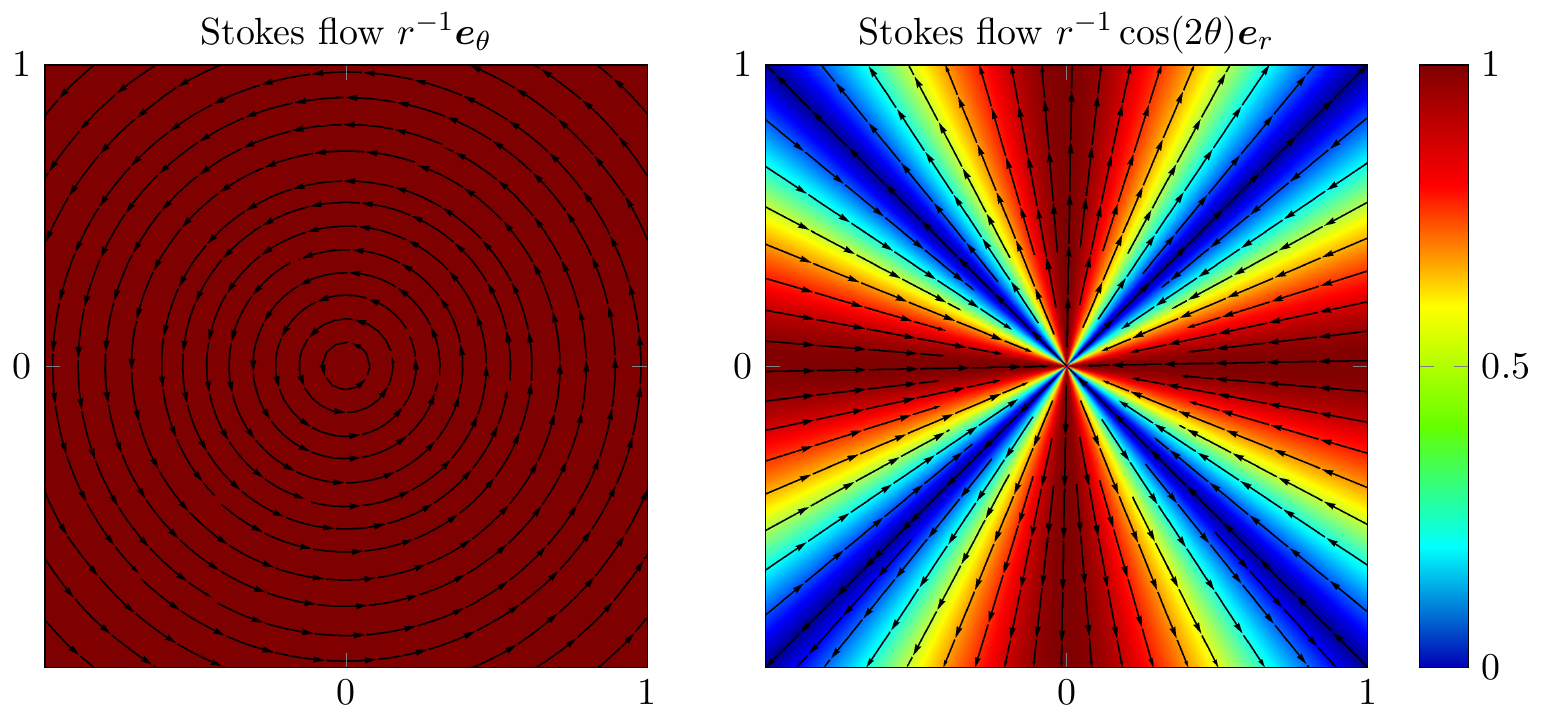}

\caption{\label{fig:stokes}Representation of the velocity vector field produced
by the two solutions of the Stokes equations decaying like $r^{-1}$
. The first one is generated by the torque $\int_{\mathbb{R}^{2}}\left(\bx\bwedge\bff\right)$
and the second one by $\int_{\mathbb{R}^{2}}\left(x_{1}f_{1}-x_{2}f_{2}\right)$.}
\end{figure}

We now discuss the consequences of the existence of the newly found
solutions for the solutions of the Navier-Stokes equations (\ref{eq:ns-force})
in $\mathbb{R}^{2}$. For the reasons explained in the introduction,
if $\bff$ has non-zero mean, the solution can not decay like $r^{-1}$.
Therefore, the new solutions describe, at best, the asymptotes of
solutions for the case where $\bff$ has zero mean. To discuss this
question, we consider the Stokes approximation,
\[
\Delta\bu-\bnabla p=\bff\,,
\]
with $\bff$ having zero mean. The asymptotic behavior of the solution
to this equation is
\[
\bu=\frac{-1}{4\pi r}\left[\left(\cos(2\theta)\int_{\mathbb{R}^{2}}\left(x_{1}f_{1}-x_{2}f_{2}\right)+\sin(2\theta)\int_{\mathbb{R}^{2}}\left(x_{1}f_{2}+x_{2}f_{1}\right)\right)\be_{r}+\int_{\mathbb{R}^{2}}\left(\bx\bwedge\bff\right)\be_{\theta}\right]+O(r^{-2})\,,
\]
By rotating the coordinates system, we can always make $\int_{\mathbb{R}^{2}}\left(x_{1}f_{2}+x_{2}f_{1}\right)=0$,
so that the Stokes approximation has two compatibility conditions:
the torque $\int_{\mathbb{R}^{2}}\left(\bx\bwedge\bff\right)$ and
$\int_{\mathbb{R}^{2}}\left(x_{1}f_{1}-x_{2}f_{2}\right)$ which are
represented in figure~\ref{fig:stokes}. Even if the small solution
found for $n=2$ has the appropriate form $\cos(2z)$ which is similar
to the $\cos(2\theta)$ of the Stokes solution, the asymptotic behavior
is likely not described by this exact solution alone, because only
one of the compatibility conditions can be lifted by this exact solution.
In addition, we note that there are two exact solutions decaying like
$r^{-1}$ and having arbitrary small torque: the harmonic function
$\mu\be_{\theta}/r$ and the new solution for $n=2$. This emphasizes
the wide variety of asymptotic behavior of the solutions to (\ref{eq:ns-force})
with small data, since by truncating one of these solutions near the
origin we obtain an exact solution in $\mathbb{R}^{2}$ with a certain
small source term. In fact, numerical studies make us believe that
even if the source term has zero mean, the solution of (\ref{eq:ns-force})
is in general not bounded by $r^{-1}$.

\phantomsection\addcontentsline{toc}{section}{\refname}\bibliography{}

\end{document}